\documentclass[12pt,draft]{amsart}
\usepackage{amsmath,amsthm,latexsym,amscd,amsbsy,amssymb,graphics}
 \relax

\chardef\bslash=`\\ 

\makeatletter
\def\verbatim{\interlinepenalty\@M \@verbatim
  \leftskip\@totalleftmargin\advance\leftskip2pc
  \frenchspacing\@vobeyspaces \@xverbatim}
\makeatother
\makeatletter
  \def\dgt@k{\dg@DX=1 \dg@DY=-4 \dg@SIZE=1}
\makeatother

\theoremstyle{plain}
\newtheorem{thm}{Theorem}[section]
\newtheorem{cor}[thm]{Corollary}
\newtheorem{lem}[thm]{Lemma}
\newtheorem{pro}[thm]{Proposition}

\theoremstyle{definition}
\newtheorem{rem}[thm]{Remark}
\newtheorem{defin}[thm]{Definition}

\newcounter{rmnum}

\newcommand{\ind}{\operatorname{ind}_{L}}
\newcommand{\Ind}{\operatorname{Ind}_{L}}

\numberwithin{equation}{section}

\begin{document}


\title[Extraordinary dimension theories generated by complexes]
{Extraordinary dimension theories generated by complexes}
\author{A.~Chigogidze}
\address{Department of Mathematics and Statistics,
University of Saskatche\-wan,
McLean Hall, 106 Wiggins Road, Saskatoon, SK, S7N 5E6,
Canada}
\email{chigogid@math.usask.ca}
\thanks{Author was partially supported by NSERC research grant.}



\begin{abstract}{We study the extraordinary dimension function $\dim_{L}$ introduced by \v{S}\v{c}epin. An axiomatic characterization of this dimension function is obtained. We also introduce inductive dimensions $\ind$ and $\Ind$ and prove that for separable metrizable spaces all three coincide. Several results such as characterization of $\dim_{L}$ in terms of partitions and in terms of mappings into $n$-dimensional cubes are presented. We also prove the converse of the Dranishnikov-Uspenskij theorem on dimension-raising maps.}
\end{abstract}

\maketitle \markboth{A.~Chigogidze}{Extraordinary dimension theories}

\section{Introduction}\label{S:intro}
In recent years there has been a significant development in Extension Theory. Number of results of classical dimension theory has been reexamined, better understood and their far reaching generalizations found.
The fundamental problem, studied in this theory, is the possibility of extending a map $f \colon A \to L$, defined on a closed subset $A$ of a space $X$, with values lying in a complex $L$, over the whole $X$ (when all such extension problems are solvable for a given space $X$ we write $L \in {\rm AE}(X)$; see Section \ref{S:pre} for precise definitions). Below we study the dimension function $\dim_{L}$ generated by a complex $L$. The dimension $\dim_{L}X$, introduced by \v{S}\v{c}epin \cite[p.984]{scepin}, is defined as the smallest integer $n$ such that $\Sigma^{n}L \in {\rm AE}(X)$, where $\Sigma^{n}L$ denotes the $n$-th iterated suspension of $L$. It is noted in \cite{scepin} that $\dim_{L}$ satisfies all of the Alexandrov's axioms except the normalization axiom (see Section \ref{S:axiomatic} for details) and by analogy with homology and cohomology theories is referred as an extraordinary dimension function with classifying complex $L$. We also introduce the small and large inductive dimensions $\ind$ and $\Ind$ generated by a complex $L$ and prove that $\dim_{L}X = \ind X = \Ind X$ for any separable metrizable space $X$ (Theorem \ref{T:coincide}). This allows us to study properties of the dimension $\dim_{L}$ by using the standard inductive approach. We would like to mention the following three characterizations.

{\bf Characterization of $\dim_{L}$ in terms of partitions (Theorem \ref{T:partition}).} {\em Let $X$ be a separable metrizable space and $n \geq 1$. Then the following conditions are equivalent:
\begin{itemize}
\item[(a)]
$\operatorname{dim}_{L}X \leq n$;
\item[(b)]
for every collection $(A_{1},B_{1}), (A_{2},B_{2}),\dots ,(A_{n},B_{n})$ of $n$ pairs of disjoint closed subsets of $X$ there exist closed sets $C_{1}, C_{2},\dots ,C_{n}$ such that $C_{i}$ is a partition between $A_{i}$ and $B_{i}$ and $L \in {\rm AE}\left(\bigcap_{i=1}^{n}C_{i}\right)$.
\end{itemize}}

 \noindent This result seems to be providing a new insight even for the standard covering dimension $\dim$. 

Another characterization of the dimension $\dim_{L}$ is contained in the following result.

{\bf Characterization of $\dim_{L}$ in terms of mappings into cubes (Theorem \ref{T:charact}).} {\em Let $X$ be a compact metrizable space and $n \geq 1$. Then the following conditions are equivalent:
\begin{itemize}
\item[(i)]
$\dim_{L}X \leq n$.
\item[(ii)]
The set of maps $f \colon X \to I^{n}$ with $\dim_{L} f =0$ forms a dense $G_{\delta}$-subset of the space $C(X,I^{n})$.
\item[(iii)]
There exists a map $f \colon X \to I^{n}$ such that $\dim_{L} f = 0$.
\end{itemize}
}

The characterizing property contained in this theorem becomes an axiom in an axiomatic characterization of $\dim_{L}$.

{\bf Axiomatic characterization of $\dim_{L}$ (Theorem \ref{T:axiomatics}).}
{\em The dimension $\dim_{L}$ is the only function,
defined on the class of finite-dimensional (in the sense of $\dim_{L}$ compact metrizable spaces, which satisfies the following axioms:
\begin{description}
\item[(C1)-- normalization axiom]
$d(X) \in \{ 0,1,2,\dots \}$ and $d(X) = 0$ if and only if $L \in {\rm AE}(X)$. 
\item[(C2) -- monotonicity axiom]
If $A$ is a closed subspace of $X$, then $d(A) \leq d(X)$.
\item[(C3) -- Poincar\'{e}'s axiom]
If $d(X) > 0$, then there exists a closed subspace $A$ in $X$ separating $X$ and such that $d(A) < d(X)$.
\item[(C4) -- Hurewicz's axiom]
If there exists a map $f \colon X \to I^{n}$ such that $d(f^{-1}(y)) =0$ for every $y \in f(X)$, then $d(X) \leq n$.
\end{description} 
}

It is important to emphasize the significance of the join and the smash product constructions in the Extension theory. Generalizations of the addition and product theorems (see Theorems \ref{T:menger} and \ref{T:product1} respectively), as well as the analog of the Hurewicz's theorem on dimension-lowering maps (Theorem \ref{T:hur}), not only well demonstrate this point, but, in fact, uncover much deeper roots of the classical prototypes of the mentioned results. In light of this it is interesting to note that within the theory of the dimension function $\dim_{L}$ the role of the smash product construction becomes superior and allows us to state the corresponding results in a much more familiar manner. Here is the illustration:

\begin{description}
\item[The Addition Theorem (Proposition \ref{P:mu})]
\[\dim_{L\wedge L}(X\cup Y) \leq \dim_{L}X +\dim_{L}Y +1 .\]
\item[The Product Theorem (Proposition \ref{P:product})]
\[\dim_{L\wedge L}(X\times Y) \leq \dim_{L}X +\dim_{L}Y .\]
\item[The Hurewicz's Inequality (Proposition \ref{P:hur})]

\[\dim_{L\wedge L}X \leq \dim_{L}Y + \dim_{L}f .\]
\end{description}

Note that the smash product $L\wedge L$ cannot be replaced by $L$ itself in neither of the above results unless $[L \wedge L] = [L]$. Of course, this is the case for $L =S^{0}$.

\section{Preliminaries}\label{S:pre}
All spaces considered below are assumed to be completely regular and Hausdorff. Letters $L$ and $K$ are reserved exclusively for locally finite countable simplicial complexes (alternatively, the reader may assume, in a majority of instances, that spaces denoted by $L$ and $K$ are Polish ${\rm ANR}$-spaces).

For a normal space $X$, the notation $L \in {\rm AE}(X)$ means that every map $f \colon A \to L$, defined on a closed subspace $A$ of $X$, admits an extension $\tilde{f} \colon X \to L$ over $X$. For a non-normal space $X$, the relation $L \in {\rm AE}(X)$ is understood in a slightly adjusted manner (see \cite{chi3}, \cite{chihand} for details). For normal spaces the modified definition coincides with the one presented above.

Following \cite{dr93}, we say that $L \leq K$ if for each space $X$ the condition $L \in {\rm AE}(X)$ implies the condition $K \in {\rm AE}(X)$. Equivalence classes of complexes (Polish ${\rm ANR}$-spaces) with respect to this relation are called extension types. The above defined relation $\leq$ creates a partial order in the class of extension types. This partial order is denoted by $\leq$ and
the extension type with representative $L$ is denoted by $[L]$. Note that under these definitions the class of all extension types has both maximal and minimal elements.
The minimal element is the extension type of the $0$-dimensional sphere $S^{0}$ and the maximal element is obviously the extension type of the one-point space $\{ {\rm pt} \}$ (or, equivalently, of any contractible complex).

$\operatorname{Cone}(L)$ and $\Sigma^{n}L$ denote respectively the cone and the $n$-th iterated suspension of $L$. $L \ast K$ and $L \wedge L$ denote the join and the smash product of $L$ and $K$.

For the reader's convenience in this section we present some facts which are needed below.

\begin{thm}[\cite{dra1}, Theorem 4]\label{T:suspension}
If $L \in {\rm AE}(X)$, then $\Sigma L \in {\rm AE}(X\times I)$.
\end{thm}

\begin{thm}[\cite{dra2}, Proposition 2.3]\label{T:countable}
If a normal space $X$ is represented as the union $X = \cup_{i=1}^{\infty}F_{i}$ of its closed subsets $F_{i}$ such that $L \in {\rm AE}(F_{i})$ for each $i$, then $L \in {\rm AE}(X)$.
\end{thm}

\begin{thm}[\cite{DD}, Theorem 5.6, Corollary 5.7]\label{T:product1}
Let $X$ and $Y$ be finite dimensional compact spaces. If $L \in {\rm AE}(X)$ and $K \in {\rm AE}(Y)$, then $L \wedge K \in {\rm AE}(X\times Y)$.
If both $L$ and $K$ are connected and $L$, in addition, is finitely dominated, then the above conclusion remains valid for all finite dimensional separable metrizable spaces.
\end{thm}

\begin{thm}[\cite{dydak1}, Theorem 1.2]\label{T:menger}
If a metrizable space $X$ is the union of two subsets  $A$, $B$ such that $L \in {\rm AE}(A)$ and $K \in {\rm AE}(B)$, then $L\ast K \in {\rm AE}(X)$.
\end{thm}

\begin{thm}[\cite{dra2}, Theorem 1; \cite{DD}, Theorem 3.7]\label{T:eilenberg}
If $L \ast K \in {\rm AE}(X)$, then for any map $f \colon A \to L$ of a closed subset $A \subseteq X$, defined on a separable metrizable space $X$, there exists a closed subset $B \subseteq X$ such that $K \in {\rm AE}(B)$ and $f$ admits an extension $\widetilde{f} \colon X \setminus B \to L$. 
\end{thm}

\begin{thm}[\cite{dra2}, Corollary 2; \cite{DD}, Theorem 3.8]\label{T:splitting}
If $X$ is a separable metrizable space such that $L \ast K \in {\rm AE}(X)$, then there exists a subset $A \subseteq X$ such that $L \in {\rm AE}(A)$ and $K \in {\rm AE}(X\setminus A)$.
\end{thm}

\begin{thm}[\cite{ols}]\label{T:ols}
For every separable metrizable space $X$ with $L \in {\rm AE}(X)$, there exists a completion $\widetilde{X}$ of $X$ such that $L \in {\rm AE}(\widetilde{X})$.
\end{thm}

\begin{thm}[\cite{chi3}, Theorem 3.5]\label{T:hewitt}
$L \in {\rm AE}(X)$ if and only if $L \in {\rm AE}(\upsilon X)$, where $\upsilon X$ denotes the Hewitt realcompactification of $X$.
\end{thm}

\begin{thm}[\cite{chi2}, Corollary 2.2]\label{T:stone}
If $L$ is a finitely dominated complex, then $L \in {\rm AE}(X)$ if and only if $L \in {\rm AE}(\beta X)$, where $\beta X$ denotes the Stone-\v{C}ech compactification of $X$.
\end{thm}

\begin{thm}[\cite{levlew}, Theorem 1.6; \cite{drus}, Theorem1.2]\label{T:hur}
Let $f \colon X \to Y$ be a map of compact spaces with $X$ finite dimensional. If $L \in {\rm AE}(Y)$ and $K \in {\rm AE}(f^{-1}(y))$ for each $y \in Y$, then $L \wedge K \in {\rm AE}(X)$.
\end{thm}

The following statement is closely related to the previous theorem. In it $X$ is not assumed to be finite dimensional.

\begin{thm}[\cite{chival}, Corollary 2.7]\label{T:chihur}
Let $f \colon X \to Y$ be a map of compact spaces. If $\dim Y \leq n$ and $L \in {\rm AE}(f^{-1}(y))$ for every $y \in Y$, then $\Sigma^{n}L \in {\rm AE}(X)$.
\end{thm}

\begin{thm}[\cite{drus}, Theorem 1.6]\label{T:raising}
Let $f \colon X \to Y$ be an onto map between metrizable compact spaces. If $L \in {\rm AE}(X)$ and $|f^{-1}(y)| \leq n+1$ for each $y \in Y$, then $\Sigma^{n}L \in {\rm AE}(Y)$.
\end{thm}

Some other statements, which also are needed below and proofs of which require spectral techniques, are included in Section \ref{S:app}.

\section{Extraordinary inductive dimensions $\ind$ and $\Ind$}\label{S:definitions}

Recall that $\Sigma^{n}L$ denotes the iterated suspension of $L$. Also, for notational convenience, we let $\Sigma^{0}L = L$. The relation $\Sigma^{n}L \in {\rm AE}(X)$ would be rewritten as $\dim_{L}X \leq n$, $n = 0,1,2,\dots$. In other words
\[ \dim_{L}X = \min\{ n \in{\mathbb N} \cup \{ 0\} \colon \Sigma^{n}L \in {\rm AE}(X)\} .\]

The idea of defining inductive dimensions with respect to classes of spaces is not new. It has been developed in \cite{aarts}, \cite{an} (see also \cite{chif}) and consists of replacing the empty set in the definition of standard inductive dimensions by elements of a given class of spaces. Our approach is very similar -- except the definition starts in the dimension zero.

\begin{defin}\label{D:large}
Let $L$ be a CW-complex and $X$ be a space. We say that 
\begin{itemize}
\item[(i)]
$\Ind X \leq 0$ if and only if $L \in {\rm AE}(X)$;
\item[(ii)]
$\Ind X \leq n$, $n \in {\mathbb N}$, if for every closed set $A \subseteq X$ and every open neighbourhood $V$ of $A$, there exists an open set $U \subseteq X$ such that $A \subseteq U \subseteq V$ and $\Ind \operatorname{Bd}U \leq n-1$;
\item[(iii)]
$\Ind X = n$ if $\Ind X \leq n$ and $\Ind X > n-1$;
\item[(iv)]
$\Ind X = \infty$ if $\Ind X > n$ for each $n = 0,1,2,\dots$.
\end{itemize}
\end{defin}

If the set $A$ in the above definition is assumed to be a singleton then we obtain the definition of the small inductive dimension $\operatorname{ind}_{L}X$. 

Note that $\operatorname{ind}_{S^{0}}X =     \operatorname{ind}X$ and $\operatorname{Ind}_{S^{0}}X = \operatorname{Ind}X$ for any space $X$. It is also clear that if $[L] \leq [K]$, then $\operatorname{ind}_{K}X \leq \operatorname{ind}_{L}X$ and $\operatorname{Ind}_{K}X \leq \operatorname{Ind}_{L}X$ for any space $X$. In particular, $\operatorname{ind}X \leq \ind X$ and $\operatorname{Ind}X \leq \Ind X$.

Of course, these definitions can be extended to higher ordinal numbers. We intend to investigate transfinite inductive dimensions $\ind$ and $\Ind$ and associated with them various types of infinite dimensional spaces (in the sense of the dimension function $\dim_{L}$) in a separate note. 


\subsection{General observations}\label{SS:general}
We record the following results for the future references.
\begin{pro}\label{P:mono}
Let $Y$ be a subspace of a space $X$. Then $\ind Y \leq \ind X$ and $\Ind Y \leq \Ind Y$ provided:
\begin{itemize}
\item[(a)]
$Y$ is an $F_{\sigma}$-subset and $X$ is normal;
\item[(b)]
$Y$ is an arbitrary subset and $X$ is perfectly normal.
\end{itemize}
\end{pro}
\begin{proof}
We proceed by induction on $\ind X$. If $\ind X=0$, then, according to Definition \ref{D:large}, $L \in {\rm AE}(X)$. By Theorem \ref{T:monotone}, $L \in {\rm AE}(Y)$, which, in turn, means that $\ind Y = 0$. Suppose that our statement is valid for spaces $X$ with $\ind X \leq n$ and consider a space with $\ind X = n+1$. For a point $y \in Y$ and its open neighbourhood $U \subseteq X$ choose a smaller neighbourhood (in $X$) $V$ such that $y \in V \subseteq \operatorname{Cl}V \subseteq U$ and $\ind \operatorname{Bd}V \leq n$. Clearly $y \in V \cap Y \subseteq \operatorname{Cl}_{Y}(V \cap Y) \subseteq U\cap Y$. Note that $\operatorname{Bd}_{Y}(V\cap Y) \subseteq \operatorname{Bd}_{X}V$. Since $\operatorname{Bd}_{Y}(V\cap Y)$ is an $F_{\sigma}$-subset of $\operatorname{Bd}_{X}V$ we may use the inductive assumption (in case (a)) and conclude that $\ind \operatorname{Bd}_{Y}(V\cap Y) \leq n$.
\end{proof}

\begin{lem}\label{L:separator}
Let $n \geq 1$ and $A$ and $B$ be disjoint closed subsets of a Lindel\"{o}f space $X$ with $\ind X \leq n$. Then there exist disjoint open subsets $G_{A}$ and $G_{B}$ in $X$ so that $A \subseteq G_{A}$, $B \subseteq G_{B}$ and $X \setminus (G_{A} \cup G_{B})$ is contained in the union $\cup_{k=1}^{\infty}F_{k}$ of closed sets $F_{k}$ with $\ind F_{k} \leq n-1$ for each $k=1,2,\dots$.
\end{lem}
\begin{proof}
Choose an open neighbourhood $U_{A}$ of $A$ such that $A \subseteq U_{A} \subseteq \operatorname{Cl}(U_{A})  \subseteq X \setminus B$.

For a point $x \in \operatorname{Cl}(U_{A})$ let $Ux$ denote an open neighbourhood of $x$ such that  $\ind\operatorname{Bd}(Ux) \leq n-1$ and $\operatorname{Cl}(Ux) \subseteq X\setminus B$. If $x \in X \setminus \operatorname{Cl}(U_{A})$ let $Ux$ denote an open neighbourhood of $x$ such that  $\ind\operatorname{Bd}(Ux) \leq n-1$ and $\operatorname{Cl}(Ux) \subseteq X \setminus \operatorname{Cl}(U_{A})$. 
Since $X$ is Lindel\"{o}f space the open cover $\{ Ux \colon x \in X\}$ contains a countable subcover $\{ U_{k} \}_{k=1}^{\infty}$.

Let $F_{k} = \operatorname{Bd}(U_{k})$. By construction, $\ind F_{k} \leq n-1$ for each $k=1,2,\dots$.

Let also 
\[ G_{A} = \bigcup\{ U_{k} \colon U_{k} \cap \operatorname{Cl}(U_{A}) \neq \emptyset\}\;\;\text{and}\;\; G_{B} = \bigcup\{ U_{k} \colon U_{k} \cap \operatorname{Cl}(U_{A}) = \emptyset\} .\]

It is not hard to verify that $A \subseteq G_{A} \subseteq U_{A}$, $B \subseteq G_{B} \subseteq X \setminus \operatorname{Cl}(U_{A})$ and $X \setminus (G_{A} \cup G_{B}) \subseteq \cup_{k=1}^{\infty}F_{k}$.
\end{proof}

\begin{pro}\label{P:lindelof}
Let $X$ be a Lindel\"{o}f space. Then $\dim_{L}X \leq \ind X$.
\end{pro}
\begin{proof}
We proceed by induction. If $\ind X = 0$, then, according to our definitions, $\dim_{L}X = 0$. Assume now that the statement holds for Lindel\"{o}f spaces with $\ind \leq n-1$ and consider a Lindel\"{o}f space $X$ such that $\ind X \leq n$.

Let $f \colon Y \to \Sigma^{n}L$ be a map defined on a closed subspace $Y$ of the space $X$. Fix an open neighbourhood $O$ of $Y$ such that $f$ is extendable over $\operatorname{Cl}O$ and denote such an extension by the same letter $f$. Let $Z$ be a functionally closed subset of $X$ such that $Y \subseteq Z \subseteq O$. Represent $\Sigma^{n}L = \Sigma (\Sigma^{n-1}L)$ as the union
of the two ``semispeheres" $L_{-}$ and $L_{+}$, each of which is a copy of $\operatorname{Con}(\Sigma^{n-1}L)$ and whose intersection $L_{0} = L_{-}\cap L_{+}$ is a copy of $\Sigma^{n-1}L$. Let $Z_{-} = f^{-1}(L_{-}) \cap Z$, $Z_{+} = f^{-1}(L_{+}) \cap Z$ and $Z_{0} = f^{-1}(L_{0}) \cap Z = Z_{-}\cap Z_{+}$. Note that $Z_{-}$, $Z_{+}$ and $Z_{0}$ also are functionally closed subsets of $X$. Consequently, $X \setminus Z_{0}$, being functionally open (and hence $F_{\sigma}$-) in $X$, is a Lindel\"{o}f space. 

Note that $Z_{-}\setminus Z_{0}$ and $Z_{+}\setminus Z_{0}$ are disjoint (functionally) closed subsets of $X\setminus Z_{0}$. By Proposition \ref{P:mono}, $\ind (X\setminus Z_{0}) \leq \ind X \leq n$. By Lemma \ref{L:separator}, there exist disjoint open subsets $G_{-}$ and $G_{+}$ in $X\setminus Z_{0}$ so that $Z_{-}\setminus Z_{0} \subseteq G_{-}$, $Z_{+}\setminus Z_{0} \subseteq G_{+}$ and $(X \setminus Z_{0}) \setminus (G_{A} \cup G_{B})$ is contained in the union $\cup_{k=1}^{\infty}F_{k}$ of closed (in $X\setminus Z_{0}$) sets $F_{k}$ with $\ind F_{k} \leq n-1$ for each $k=1,2,\dots$. By the inductive assumption, $\dim_{L}F_{k} \leq n-1$ for each $k=1,2,\dots$.

Next consider the following closed subsets of $X$: $X_{-} = (X \setminus G_{+}) \cup Z_{0}$, $X_{+} = (X \setminus G_{-}) \cup Z_{0}$ and $X_{0} = X_{-} \cap X_{+}$. Observe that $X_{-} \cap Z = Z_{-}$, $X_{+} \cap Z = Z_{+}$, $X_{0} \cap Z = Z_{0}$. Note also that $X_{0} \setminus Z_{0} \subseteq \cup_{k=1}^{\infty}F_{k}$. Theorem \ref{T:countable}
 guarantees that $\dim_{L}(X_{0} \setminus Z_{0}) \leq \dim_{L}\left(\cup_{k=1}^{\infty}F_{k}\right) \leq n-1$. The map $f|Z_{0} \colon Z_{0} \to L_{0} = \Sigma^{n-1}L$ admits an extension
$g_{0} \colon \operatorname{Cl}G \to L_{0}$, where $G$ is an open neighbourhood of $Z_{0}$ in $X_{0}$. Since, as was noted, $\dim_{L}(X_{0}\setminus Z_{0}) \leq n-1$, it follows that the maps 

\[ g_{0}|\left[ ( X_{0} \setminus Z_{0}) \cap \operatorname{Cl}G \right] \colon (X_{0} \setminus Z_{0}) \cap \operatorname{Cl}G \to L_{0}\]

\noindent  admits an extension $g \colon X_{0} \setminus Z_{0} \to L_{0}$ onto the whole  $X_{0} \setminus Z_{0}$. Let $h_{0} \colon X_{0} \to L_{0}$ be the map which coincides with $f$ on $\operatorname{Cl}G$ and with $g_{0}$ on $X_{0} \setminus G$.

Next consider the map $h_{-} \colon X_{-} \to L_{-}$, defined by letting
\[ h_{-}(x) = 
\begin{cases}
f(x) , \text{if}\; x \in Z_{-}\\
h_{0}(x) , \text{if}\;  x \in X_{0}.\\
\end{cases}\]

Since $L_{-}$ is a contractible complex, the map $h_{-}$ can be extended to a map $H_{-} \colon X_{-} \to L_{-}$.

Similarly the map $h_{+} \colon X_{+} \to L_{+}$, defined by letting
\[ h_{+}(x) = 
\begin{cases}
f(x) , \text{if}\; x \in Z_{+}\\
h_{0}(x) , \text{if}\;  x \in X_{0},\\
\end{cases}\]

\noindent also admits an extension $H_{+} \colon X_{+} \to L_{+}$. Note that the maps $H_{-}$ and $H_{+}$ agree on $X_{0}$ (with the map $h_{0}$) and hence define the map $H \colon X \to \Sigma^{n}L$, which obviously is an extension of the originally given map $f$. This proves that $\dim_{L}X \leq n$.
\end{proof}

\begin{rem}\label{R:different}
Let $L$ be a connected non-contractible CW complex such that the extension type $[L]$ is bounded from above by the extension type of some sphere.  Without loss of generality we may assume that $[L] \neq [S^{m}]$ for any $m$. Let $n$ be the smallest integer such that $[L] < [S^{n}]$. We show that there exists a compact space $X$ such that $\dim_{L}X$ is finite, but $\ind X$ is not even defined. Indeed, consider an $(n+2)$-dimensional compact space $X_{n}$ without intermediate dimensions (i.e. if $F$ is a closed subspace of $X_{n}$, then either $\dim F = n+2$ or $\dim F = 0$). Let $k$ be the smallest non-negative integer such that the homotopy group $\pi_{k+1}(L)$ is non-trivial. Note that then $[S^{k}] \leq [L]$ and consequently $[S^{n+2}] = [\Sigma^{n+2-k}\left( S^{k}\right)] \leq [\Sigma^{n+2-k}L]$. This implies that $\dim_{L}X_{n} \leq n+2-k$. Note also that $\dim_{L}X_{n} \geq 1$ (to see this observe that $\dim_{L}X_{n} = 0$ simply means that $L \in {\rm AE}(X)$ which, in light of $[L] < [S^{n}]$, would imply $\dim X \leq n$ contradicting the choice of the compactum $X_{n}$).
Next suppose that the small inductive dimension $\ind X_{n}$ is finite, i.e. $p = \ind X_{n} \geq \dim_{L}X_{n} >0$. Then $X_{n}$ contains a closed subset $F$ such that $\ind F = 1$. By Proposition \ref{P:lindelof}, $\dim_{L}F \leq \ind F = 1$. Since $[\Sigma L] \leq S^{n+1}$ we conclude that $\dim F \leq n+1$. But $X_{n}$ does not contain positive-dimensional closed subsets of the covering dimension strictly less than $n+2$. Consequently $\dim F = 0$, which is impossible in view of $\ind F = 1$. 
\end{rem}

\begin{pro}\label{P:stone}
Let $L$ be a finitely dominated. Then $\Ind X = \Ind\beta X$ for any normal space $X$.
\end{pro}
\begin{proof}
First we show that $\Ind X \leq \Ind\beta X$. If $\Ind\beta X = 0$, then $L \in {\rm AE}(\beta X)$ which, by Theorem \ref{T:stone}, implies that $L \in {\rm AE}(X)$. Consequently, $\Ind X =0$. Suppose that the inequality $\Ind X \leq \Ind\beta X$ is valid for all normal spaces $X$ with $\Ind\beta X \leq n-1$ and consider a normal space $X$ such that $\Ind\beta X = n$. Let $A$ and $B$ are disjoint closed subsets of $X$. Then $\operatorname{Cl}_{\beta X}A \cap \operatorname{Cl}_{\beta X}B = \emptyset$. Choose open subsets $U$ and $V$ in $\beta X$ so that $\operatorname{Cl}_{\beta X}A \subseteq U$, $\operatorname{Cl}_{\beta X}B \subseteq V$ and $\Ind F \leq n-1$ where $F = \beta X \setminus (U \cup V)$. Clearly the set $X \cap F$ separates the sets $A$ and $B$ in $X$. Note also that $X \cap F$ is normal. Since $\beta F \cap X = \operatorname{Cl}_{\beta X}(F\cap X) \subseteq F$ it follows from the inductive assumption that $\Ind (F\cap X) \leq n-1$. This proves that $\Ind X \leq n$.

Next we prove the inequality $\Ind\beta X \leq \Ind X$. If $\Ind X = 0$, then $L \in {\rm AR}(X)$ and by Theorem \ref{T:stone}, $L \in {\rm AE}(\beta X)$. This means that $\Ind \beta X = 0$. Assume that the inequality $\Ind\beta X \leq \Ind X$ holds for all normal spaces with $\Ind X \leq n-1$ and consider a normal space $X$ with $\Ind X = n$. Let $A$ and $B$ be disjoint closed subsets of $\beta X$. Consider open subsets $U$ and $V$ in $\beta X$ such that $A \subseteq U \subseteq \operatorname{Cl}U$, $B \subseteq V \subseteq \operatorname{Cl}V$ and $\operatorname{Cl}U \cap \operatorname{Cl}V =\emptyset$. Then $X \cap \operatorname{Cl}U$ and $X \cap \operatorname{Cl}V$ are disjoint nonempty closed subsets of $X$. Since $\Ind X \leq n$, these closed sets can be separated (in $X$) by a closed subset $F \subseteq X$ with $\Ind F \leq n-1$. Obviously $\operatorname{Cl}_{\beta X}F$ separates $\operatorname{Cl}_{\beta X}(X \cap \operatorname{Cl}_{\beta X}U)$ and $\operatorname{Cl}_{\beta X}(X \cap \operatorname{Cl}_{\beta X}V)$. Note also that
\[ A \subseteq \operatorname{Cl}_{\beta X}U \subseteq \operatorname{Cl}_{\beta X}(X \cap U) \subseteq \operatorname{Cl}_{\beta X}(X \cap \operatorname{Cl}_{\beta X}U)\]

\noindent and

\[ B \subseteq \operatorname{Cl}_{\beta X}V \subseteq \operatorname{Cl}_{\beta X}(X \cap V) \subseteq \operatorname{Cl}_{\beta X}(X \cap \operatorname{Cl}_{\beta X}V) .\]

This shows that $\operatorname{Cl}_{\beta X}F$ is a separator between the sets $A$ and $B$ in $\beta X$. Since $F$ is normal and since $\beta F = \operatorname{Cl}_{\beta X}F$ the inductive assumption implies that $\Ind \operatorname{Cl}_{\beta X}F \leq \Ind F \leq n-1$. This shows that $\Ind \beta X \leq n$.
\end{proof}

\begin{cor}\label{C:vedenosov}
Let $L$ be a finitely dominated. If $X$ is a normal space, then $\dim_{L}X \leq \Ind X$.
\end{cor}
\begin{proof}
By Proposition \ref{P:stone}, $\ind X = \ind\beta X$. Since $L$ is finitely dominated, so are its iterated suspensions and consequently, by \cite[Corollary 2.2]{chi2}, $\dim_{L}X = \dim_{L}\beta X$ . Then, according to Proposition \ref{P:lindelof}, we have
\[ \dim_{L}X = \dim_{L}\beta X \leq \ind\beta X \leq \Ind\beta X = \Ind X .\]
\end{proof}


\subsection{Inductive dimensions $\ind$ and $\Ind$ of perfectly normal spaces}\label{SS:perfect}

Basic properties of classical inductive dimensions $\operatorname{ind}$ and $\operatorname{Ind}$ have their prototypes for the dimensions $\ind$ and $\Ind$ in perfectly normal spaces. Proofs of the following two statements are standard and require only straightforward adjustments based on Theorems \ref{T:countable} and \ref{T:monotone}(b).

\begin{thm}\label{T:indcountable}
If a perfectly normal space $X$ can be represented as the union of a countable collection $X = \bigcup_{i=1}^{\infty}X_{i}$ of closed subsets such that $\Ind X_{i} \leq n$ for each $i = 1,2,\dots$, then $\Ind X \leq n$.
\end{thm}

\begin{thm}\label{T:Indmonotone}
If $Y$ is a subspace of a perfectly normal space $X$, then $\Ind Y \leq \Ind X$.
\end{thm}

Theorems \ref{T:indcountable} and \ref{T:Indmonotone} have several corollaries. Some of them, proofs of which follow standard schemes, are presented below.

\begin{pro}\label{P:indstone}
Let $L$ be a finitely dominated complex. If $X$ is a perfectly normal space, then $\ind \beta X = \Ind X = \Ind \beta X$.
\end{pro}

\begin{pro}\label{P:lindperfect}
If $X$ is a perfectly normal Lindel\"{o}f space, then $\ind X = \Ind X$.
\end{pro}

\begin{pro}[cf. Proposition \ref{P:mu}]\label{T:mengerurysohn}
If $X \cup Y$ is perfectly normal, then $\operatorname{ind}_{L \ast L}(X \cup Y) \leq \ind X +\ind Y$ and $\operatorname{Ind}_{L\ast L}(X\cup Y) \leq \Ind X + \Ind Y$.
\end{pro}
\begin{proof}
Let us prove the inequality $\operatorname{Ind}_{L\ast L}(X\cup Y) \leq \Ind X + \Ind Y$. Proof of the remaining one is similar. We proceed by induction with respect to the number $n = \Ind X + \Ind Y$. First consider the case $n=0$. In this case, $L \in {\rm AE}(X)$, $L \in {\rm AE}(Y)$ and according to Theorem \ref{T:join}, $L \ast L \in {\rm AE}(X\cup Y)$. This simply means that $\operatorname{Ind}_{L\ast L}(X \cup Y) =0$. 

Next suppose that the inequality is correct in situations when $n = \Ind X + \Ind Y \leq m$ for some $m \geq 0$ and consider the case with $n = m+1$. Without loss of generality we may assume that $\Ind X \geq 1$. Let $A$ and $B$ be disjoint closed subsets in $X \cup Y$. Choose a closed subset $C \subseteq X\cup Y$ separating $A$ and $B$ and such that $\Ind (C \cap X) < \Ind X$. Then the perfectly normal space $C$ is represented as the union of two subsets $C = (C \cap X) \cup (C \cap Y)$. As noted, by the choice of the set $C$ and by Theorem \ref{T:Indmonotone}, $\Ind (C \cap X) + \Ind (C \cap Y) < \Ind X + \Ind Y = m+1$. By the inductive assumption, $\operatorname{Ind}_{L\ast L}C \leq \Ind (C \cap X) + \Ind (C \cap Y) \leq m$. This proves that $\operatorname{Ind}_{L\ast L}(X \cup Y) \leq m+1 = \Ind X + \Ind Y$.
\end{proof}


\subsection{Inductive dimensions $\ind$ and $\Ind$ of separable metrizable spaces}\label{S:separable}

The following statement expresses a basic fact connecting all three dimension functions $\dim_{L}$, $\ind$ and $\Ind$.
\begin{thm}\label{T:coincide}
If $X$ is a separable metrizable space, then $\dim_{L}X = \ind X =$ $\Ind X$.
\end{thm}
\begin{proof}
The inequality $\ind X \leq \Ind X$ trivially holds for any space and the inequality $\dim_{L}X \leq \ind X$, according to Proposition \ref{P:lindelof}, is true for all Lindel\"{o}f spaces. 

Let us prove the remaining inequality $\Ind X \leq \dim_{L}X$. We proceed by induction. Clearly $\dim_{L}X =0$ implies $\Ind X = 0$. Assume that the inequality is valid for separable metrizable spaces $Y$ with $\dim_{L}Y \leq n-1$ and consider a space $X$ with $\dim_{L}X = n$, $n \geq 1$. Let $A$ be a closed subset of $X$ and $V$ be its open neighbourhood. Consider the map $f \colon A \cup (X \setminus V) \to S^{0}$ with $f(A) = 0$ and $f(X\setminus V) = 1$. Note that $\Sigma^{n}L$ is canonically homeomorphic to the join $S^{0}\ast \Sigma^{n-1}L$. By Theorem \ref{T:eilenberg}, $f$ can be extended to a map $g \colon X\setminus Y \to S^{0}$, where $Y$ is a closed subset in $X$ such that $\dim_{L}Y \leq n-1$. By the inductive assumption, $\Ind Y \leq n-1$. Let $U = g^{-1}(0)$. Obviously $U$ is open in $X$ and $A \subseteq U \subseteq \operatorname{Cl}U \subseteq U \cup Y \subseteq X \setminus g^{-1}(1) \subseteq V$. It only remains to note that $\operatorname{Bd}U \subseteq Y$ and consequently, by Proposition \ref{P:mono}, $\Ind\operatorname{Bd}U \leq n-1$. 
\end{proof}

\begin{thm}\label{T:decomposition}
Let $X$ be a separable metrizable space and $n \geq 0$. Then the following conditions are equivalent:
\begin{itemize}
\item[(a)]
$\ind X \leq n$;
\item[(b)]
$X$ can be represented as the union $X = X_{1} \cup X_{2} \cup \dots \cup X_{n+1}$, where $L \in {\rm AE}(X_{1})$ and $\dim X_{k} \leq 0$ for each $k = 2,\dots ,n+1$.
\end{itemize}
\end{thm}
\begin{proof}
(a) $\Rightarrow$ (b). If $\ind X = 0$ the statement is trivially true. Assume that the implication is true for all separable metrizable spaces $X$ satisfying the inequality $\ind X \leq n-1$, $n \geq 1$, and consider a space $X$ such that $\ind X = n$. Take a countable open base ${\mathcal U} = \{ U_{i} \colon i \in \mathbb{N} \}$ of $X$ such that $\ind \operatorname{Bd}(U_{i}) \leq n-1$ for each $i \in \mathbb{N}$. By Theorems \ref{T:coincide} and \ref{T:countable}, $\ind\left(\cup\{ \operatorname{Bd}(U_{i}) \colon i \in \mathbb{N}\}\right) \leq n-1$. By the inductive assumption $\cup\{ \operatorname{Bd}(U_{i}) \colon i \in \mathbb{N}\} = X_{1}\cup X_{2} \cup \dots \cup X_{n}$, where $L \in {\rm AE}(X_{1})$ and $\dim X_{i} = 0$ for each $i = 2,\dots ,n$. Obviously the subspace $X \setminus \cup\{ \operatorname{Bd}(U_{i}) \colon i \in \mathbb{N}\}$, as a space with base consisting of open and closed subsets, is zero-dimensional. Then $X = X_{1} \cup X_{2} \cup \dots \cup X_{n} \cup X_{n+1}$, where $X_{n+1} = X \setminus \cup\{ \operatorname{Bd}(U_{i}) \colon i \in \mathbb{N}\}$, is the needed decomposition of $X$. 

(b) $\Rightarrow$ (a). Clearly $\dim \left( \cup_{i=2}^{n+1}X_{i}\right) \leq n-1$. In other words, $S^{n-1} \in {\rm AE}\left( \cup_{i=2}^{n+1}X_{i}\right)$. By Theorem \ref{T:menger}, $L \ast S^{n-1} \in {\rm AE}(X)$. Since $[\Sigma^{n}L] = [L\ast S^{n-1}]$, it follows that $\dim_{L}X \leq n$. Theorem \ref{T:coincide} completes the proof.
\end{proof}


\section{Further properties of the dimension function $\dim_{L}$}\label{S:further}

In this section we present several statements related to the dimension $\dim_{L}$.


\subsection{Characterization of $\dim_{L}$ in terms of partitions}\label{SS:partitions}
The following is a prototype of the classical characterization of the dimension $\dim$ in terms of partitions. 

\begin{thm}\label{T:partition}
Let $X$ be a separable metrizable space and $n \geq 1$. Then the following conditions are equivalent:
\begin{itemize}
\item[(a)]
$\operatorname{dim}_{L}X \leq n$;
\item[(b)]
for every collection $(A_{1},B_{1}), (A_{2},B_{2}),\dots ,(A_{n},B_{n})$ of $n$ pairs of disjoint closed subsets of $X$ there exist closed sets $C_{1}, C_{2},\dots ,C_{n}$ such that $C_{i}$ is a partition between $A_{i}$ and $B_{i}$ and $L \in {\rm AE}\left(\bigcap_{i=1}^{n}C_{i}\right)$.
\end{itemize}
\end{thm}
\begin{proof}
(a) $\Longrightarrow$ (b). By Theorems \ref{T:coincide} and  \ref{T:decomposition}, $X$ can be represented as the union $X = X_{1} \cup X_{2} \cup \dots \cup X_{n+1}$, where $L \in {\rm AE}(X_{1})$ and $\dim X_{k} \leq 0$ for each $k = 2,\dots ,n+1$. By \cite[Theorem 1.2.11]{eng}, for each $i = 1,\dots ,n$ there exists a partition $C_{i}$ between $A_{i}$ and $B_{i}$ such that $C_{i} \cap X_{i+1} = \emptyset$. Obviously $\bigcap_{i=1}^{n}C_{i} \subseteq X \setminus \bigcup_{i=1}^{n}X_{i+1} \subseteq X_{1}$. Since $L \in {\rm AE}(X_{1})$ it follows that $L \in {\rm AE}\left( \bigcap_{i=1}^{n}C_{i}\right)$.

(b) $\Longrightarrow$ (a). The iterated suspension $\Sigma^{n}L$ is canonically homeomorphic to the join $L \ast S^{n-1}$, which, in turn, is canonically homeomorphic to the iterated join $L \ast S^{0}_{1}\ast \cdots \ast S_{n}^{0}$, where $S_{i}^{0} = \{ s_{0}^{i},s_{1}^{i}\}$, $i = 1,\dots ,n$, is a copy of the zero-dimensional sphere $S^{0} = \{ s_{0},s_{1}\}$. This iterated join is homeomorphic to the subspace $\widetilde{L}$ of the product $\operatorname{Cone}(L)\times \prod_{i=1}^{n}\operatorname{Cone}\left( S_{i}^{0}\right)$, defined by letting (see \cite[pp.185--188]{post} and \cite[pp. 48--50]{rf} for details)
\begin{multline*}
 \widetilde{L} = \big\{ \left([l,t_{0}], \left([x_{i},t_{i}]\right)_{i=1}^{n}\right) \in\\ \operatorname{Cone}(L)\times \prod_{i=1}^{n}\operatorname{Cone}\left( S_{i}^{0}\right) \colon t_{i} = 1\; \text{for at least one}\; i = 0,1,\dots ,n \big\} .
\end{multline*}

\noindent In other words, $\widetilde{L}$ is the union of ``faces" of the entire product.

Note that there exists a retraction (``central projection") 
\[ r \colon \left(\operatorname{Cone}(L)\times \prod_{i=1}^{n}\operatorname{Cone}\left( S_{i}^{0}\right)\right) \setminus \left( v_{0}, v_{1},\dots ,v_{n}\right) \to \widetilde{L} ,\]

\noindent where $v_{0} = (L \times [0,1])/(L\times \{ 0\})$ and $v_{i} = (S_{i}^{0}\times [0,1])/(S_{i}^{0}\times \{ 0\})$, $i = 1,\dots ,n$, are the vertices of the cones $\operatorname{Cone}(L)$ and $\operatorname{Cone}\left( S_{i}^{0}\right)$ respectively.

Let also
\[ \pi_{0} \colon \operatorname{Cone}(L)\times \prod_{i=1}^{n}\operatorname{Cone}\left( S_{i}^{0}\right) \to \operatorname{Cone}(L) \]

\noindent and 
\[ \pi_{i} \colon \operatorname{Cone}(L)\times \prod_{i=1}^{n}\operatorname{Cone}\left( S_{i}^{0}\right) \to \operatorname{Cone}\left( S_{i}^{i}\right) , i = 1,\dots ,n,\]

\noindent denote the standard projections onto the corresponding coordinates.

Next consider a map $f \colon A \to \widetilde{L}$, defined on a closed subset $A$ of a space $X$. In order to prove our statement it suffices to extend $f$ over the whole $X$. Since the product $\operatorname{Cone}(L)\times \prod_{i=1}^{n}\operatorname{Cone}\left( S_{i}^{0}\right)$ is an absolute extensor there exists an extension $F \colon X \to \operatorname{Cone}(L)\times \prod_{i=1}^{n}\operatorname{Cone}\left( S_{i}^{0}\right)$  of $f$ over the whole $X$.

For each $i = 1,\dots ,n$ consider disjoint closed sets $A_{i} = F^{-1}\left( \pi_{i}^{-1}([s_{0}^{i}, 1])\right)$ and $B_{i} = F^{-1}\left(\pi_{i}^{-1}([s_{1}^{i},1])\right)$ of the space $X$. According to condition (b), for each $i = 1,\dots ,n$, there exists a closed partition $C_{i}$ in $X$ between the sets $A_{i}$ and $B_{i}$ such that $L \in {\rm AE}\left(\bigcap_{i=1}^{n}C_{i}\right)$. Choose a function $g_{i} \colon X \to \operatorname{Cone}\left( S_{i}^{0}\right)$ such that 

\[ g_{i}(x) = \begin{cases}
\left[ s_{0}^{i},1\right] ,\;\; \text{if and only if}\;\;  x \in A_{i},\\
\;\;\;\;\; v_{i}\; ,\;\;\text{if and only if}\;\; x \in C_{i} ,\\
\left[ s_{1}^{i},1\right] ,\;\; \text{if and only if}\;\;  x \in B_{i}.
\end{cases}
\]

Next consider the closed subspace $C = \bigcap_{i=1}^{n}C_{i} = \bigcap_{i=1}^{n}g_{i}^{-1}(v_{i})$ of $X$. Since $L \in {\rm AE}\left( C\right)$, it follows that the restriction 
\[ \pi_{0}\circ F|(A\cap C) \colon A \cap C \to L\times \{ 1\} \subseteq \operatorname{Cone}(L)\]

\noindent admits an extension $h \colon C \to L \times\{ 1\}$. Finally let
$g_{0} \colon X \to \operatorname{Cone}(L)$ be an extension of the map 
$\widetilde{h} \colon A \cup C \to L\times\{ 1 \} \subseteq \operatorname{Cone}(L)$, defined by letting

\[ h(x) =
\begin{cases}
\pi_{0}\circ F(x) ,\;\text{if}\; x \in A ,\\
h(x) ,\; \text{if}\; x \in C .\\
\end{cases}
\]

\noindent Note that $g_{0}|A = \pi_{0}\circ F|A$ and $g_{0}(C) \subseteq  \operatorname{Cone}(L) \setminus \{ v_{0}\}$.

The diagonal product $g(x) = \left( g_{0}(x), \dots ,g_{n}(x)\right)$, $x \in X$, defines the map $g \colon X \to \operatorname{Cone}(L) \times \prod_{i=1}^{n}\operatorname{Cone}\left( S_{i}^{0}\right)$ such that $\pi_{i}\circ g = g_{i}$ for each $i = 1,\dots ,n$. Since $(v_{i})_{i=1}^{n} \not\in g(X)$ we conclude that the composition $\widetilde{g} = r\circ g \colon X \to \widetilde{L}$ is well defined.
Note that $\widetilde{g}|A \simeq f$ as maps into $\widetilde{L}$. The corresponding homotopy can be defined by assigning to every $x \in A$ and each number $t \in [0,1]$ the point $H(x,t)$ which divides the interval (along the paths of the corresponding cones constituting the ``faces" of the product $\operatorname{Cone}(L)\times\prod_{i=1}^{n}\operatorname{Cone}\left(S_{i}^{0}\right)$ forming the subspace $\widetilde{L}$) with end-points $f(x)$ and $g(x)$ in the ratio of $t$ to $1-t$. Consequently, by the Homotopy Extension Theorem, there exists the required extension $\widetilde{f} \colon X \to \widetilde{L}$ of the originally given map $f$.
\end{proof}


\subsection{Dimensional properties of unions and products}\label{SS:union}

Extensional properties of the union $X \cup Y$ are well understood (see Theorems \ref{T:menger}, \ref{T:eilenberg} and \ref{T:splitting}). Theorem \ref{T:join} allows us to give a more familiar form to the union theorem for the dimension $\dim_{L}$.

\begin{pro}\label{P:mu}
Let $X$ and $Y$ be $z$-embedded subsets of the union $X \cup Y$. Then $\dim_{L\wedge L}(X \cup Y) \leq \dim_{L}X +\dim_{L}Y +1$.
\end{pro}
\begin{proof}
Let $\dim_{L}X = n$ and $\dim_{L}Y = m$. Then $\Sigma^{n}L \in {\rm AE}(X)$ and $\Sigma^{m}L \in {\rm AE}(Y)$. By, Theorem \ref{T:join}, $\Sigma^{n}L \ast \Sigma^{m}L \in {\rm AE}(X \cup Y)$. Next note that
\begin{multline*}
 \left[\Sigma^{n}L \ast \Sigma^{m}L \right] = \left[ \left( L \ast S^{n-1}\right) \ast \left( L \ast S^{m-1}\right) \right] = \left[ L \ast L \ast S^{n-1}\ast S^{m-1}\right] =\\ \left[ L \ast L \ast S^{n+m-1}\right] = \left[ \Sigma^{n+m}(L \ast L)\right] = \left[ \Sigma^{n+m}\left(\Sigma (L \wedge L)\right) \right] = \left[ \Sigma^{n+m+1}(L \wedge L)\right] .
\end{multline*}

Consequently, $\Sigma^{n+m+1}(L \wedge L) \in {\rm AE}(X \cup Y)$. The latter, by definition, means that 
\[ \dim_{L \wedge L}(X \cup Y) \leq n + m + 1 = \dim_{L}X +\dim_{L}Y +1 .\]
\end{proof}

\begin{cor}\label{C:mu}
Under the assumptions of Proposition \ref{P:mu}, $\dim_{L}(X\cup Y) \leq \dim_{L}X +\dim Y +1$.
\end{cor}
\begin{proof}
Note that $\dim Y = \dim_{S^{0}}Y$ and that $[L \wedge S^{0}] = [L]$.
\end{proof}

\begin{rem}
Generally speaking, the inequality $\dim_{L\wedge L}(X \cup Y) \leq \dim_{L}X +\dim_{L}Y +1$ cannot be replaced by $\dim_{L}(X \cup Y) \leq \dim_{L}X +\dim_{L}Y +1$. Indeed, let $L = S^{k}$ with $k > 0$. Let also $\dim X = n > k$, $\dim Y = m > k$ and $\dim (X\cup Y) = n+m+1$. Then 

\begin{multline*}
\dim_{L}(X \cup Y) = (n+m+1) -k > (n-k) + (m-k) +1 =\\ \dim_{L}X +\dim_{L}Y +1 .
\end{multline*}
\end{rem}

Theorem \ref{T:product1} also can be given a more familiar form.

\begin{pro}\label{P:product}
Let $L$ be finitely dominated. If $X$ and $Y$ are finitely dimensional and their product $X \times Y$ is Lindel\"{o}f, then $\dim_{L\wedge L}(X \times Y) \leq \dim_{L}X +\dim_{L}Y$.
\end{pro}
\begin{proof}
Let $\dim_{L}X = n$ and $\dim_{L}Y = m$. Recall that this means that $\Sigma^{n}L \in {\rm AE}(X)$ and $\Sigma^{m}L \in {\rm AE}(Y)$ respectively. Since $L$ is finitely dominated, we conclude, by Theorem \ref{T:stone}, that $\Sigma^{n}L \in {\rm AE}(\beta X)$ and $\Sigma^{m}L \in {\rm AE}(\beta Y)$. By Theorem \ref{T:product1}, $\left(\Sigma^{n}L\right) \wedge \left(\Sigma^{m}L\right) \in {\rm AE}(\beta X\times \beta Y)$. Finally, by Theorem \ref{T:monotone}(d), $\left(\Sigma^{n}L\right) \wedge \left(\Sigma^{m}L\right) \in {\rm AE}(X\times Y)$.
Next note that $[ \left(\Sigma^{n}L\right) \wedge \left(\Sigma^{m}L\right)] = [ \Sigma^{n+m}(L \wedge L)]$. Indeed, since $[\Sigma L] = [L \wedge S^{1}]$, we have

\begin{multline*}
 [\Sigma^{n}L \wedge \Sigma^{m}L] = [ ( L \wedge \underbrace{S^{1} \wedge \cdots \wedge S^{1}}_{n}) \wedge (L \wedge \underbrace{S^{1}\wedge \cdots \wedge S^{1}}_{m}) ] =\\ [ L\wedge L \wedge \underbrace{S^{1} \wedge \cdots \wedge S^{1}}_{n+m}] = [\Sigma^{n+m}(L \wedge L)] .
\end{multline*}

Consequently, $\Sigma^{n+m}(L \wedge L) \in {\rm AE}(X\times Y)$. This means that $\dim_{L\wedge L}(X\times Y) \leq n+m$.
\end{proof}

\begin{cor}\label{C:product}
Under the assumptions of Proposition \ref{P:product}, $\dim_{L}(X \times Y) \leq \dim_{L}X +\dim Y$.
\end{cor}
\begin{proof}
Note that $\dim Y = \dim_{S^{0}}Y$ and that $[L \wedge S^{0}] = [L]$.
\end{proof}


\subsection{Mappings and dimension}\label{SS:mappings}
Hurewicz's theorem on dimension-lowering maps also has a familiar appearance for the dimension function $\dim_{L}$. As usual, for a map $f \colon X \to Y$ we let
\[ \dim_{L} f = \sup\{ \dim_{L}f^{-1}(y) \colon y \in Y\} .\]

\begin{pro}\label{P:hur}
Let $f \colon X \to Y$ be a map of metrizable compacta with $X$ finite dimensional. Then $\dim_{L\wedge L}X \leq \dim_{L}Y +\dim_{L}f$.
\end{pro}
\begin{proof}
This is a particular case of Theorem \ref{T:hur} (alternatively, under an additional assumption of finite dimensionality of $Y$, one can use Proposition \ref{P:product} and \cite[Corollary 3.2]{brochi}). Indeed, let $\dim_{L}Y = n$ and $\dim_{L}f  = m$. Then $\Sigma^{n}L \in {\rm AE}(Y)$ and $\Sigma^{m}L \in {\rm AE}(f^{-1}(y))$ for each $y \in Y$. By the cited result, $\left(\Sigma^{n}L\right) \wedge \left(\Sigma^{m}L\right) \in {\rm AE}(X)$. As in the proof of Proposition \ref{P:product}, the latter implies the required inequality $\dim_{L\wedge L}X \leq n+m$.
\end{proof}

\begin{pro}\label{P:cubes}
Let $\dim_{L}X \leq n$. Then the set of maps $f \colon X \to I^{n}$ such that $\dim_{L}f = 0$ forms a dense $G_{\delta}$-subset in the space $C(X,I^{n})$ of all continuous maps of $X$ into the cube $I^{n}$ equipped with the compact open topology.
\end{pro}
\begin{proof} By Theorem \ref{T:decomposition}(b), $X = X_{1} \cup X_{2}$  such that $\dim_{L}X_{1} = 0$ and $\dim X_{2} \leq n-1$. By Theorem \ref{T:ols}, we may assume that $X_{1}$ is a $G_{\delta}$-subset of $X$. Then $X\setminus X_{1}$ can be written as the union of an increasing sequence $B_{1} \subseteq B_{2}\subseteq \cdots$ of closed at most $(n-1)$-dimensional subsets of $X$. By Hurewicz's theorem \cite[Ch. IV, \S 45(VIII)]{kur}, the subset

\[ C_{i} = \{ g \in C(X,I^{n}) \colon g|B_{i}\; \text{is of order}\leq n\}\]

\noindent is dense (and $G_{\delta}$) in the space $C(X,I^{n})$ for each $i$ (the order of a map $f$ does not exceed $k$ if the cardinality of each fiber is at most $k+1$). Then the intersection $C = \cap C_{i}$ is still dense (and $G_{\delta}$) in $C(X,I^{n})$. Note that for any $g \in C$, the order of the restriction $g|(X\setminus X_{1}) \colon X\setminus X_{1} \to I^{n}$ does not exceed $n$, i.e. $|g^{-1}(y) \cap (X \setminus X_{1})| \leq n$ for each $y \in I^{n}$. For any such $g$ we have $g^{-1}(y) = \left( g^{-1}(y) \cap (X \setminus X_{1})\right) \cup \left( g^{-1}(y) \cap X_{1}\right)$. Since $\dim_{L}\left( g^{-1}(y) \cap X_{1}\right) \leq \dim_{L}X_{1} = 0$, it follows that $\dim_{L}g^{-1}(y) = 0$ for any $y \in I^{n}$.
\end{proof}

The following statement provides a characterization of spaces with  dimension $\dim_{L}$ not exceeding $n$. It will be used in Section \ref{S:axiomatic}.

\begin{thm}\label{T:charact}
Let $X$ be a compact metrizable space and $n \geq 1$. Then the following conditions are equivalent:
\begin{itemize}
\item[(i)]
$\dim_{L}X \leq n$.
\item[(ii)]
The set of maps $f \colon X \to I^{n}$ with $\dim_{L} f =0$ forms a dense $G_{\delta}$-subset of the space $C(X,I^{n})$.
\item[(iii)]
There exists a map $f \colon X \to I^{n}$ such that $\dim_{L} f = 0$.
\end{itemize}
\end{thm}
\begin{proof}
The implication (i) $\Longrightarrow$ (ii) is proved in Proposition \ref{P:cubes} and the implication (ii) $\Longrightarrow$ (iii) is trivial. The remaining implication (iii) $\Longrightarrow$ (i) is contained in Theorem \ref{T:chihur}.
\end{proof}

The following statement is the converse of Dranishnikov-Uspenskij result (Theorem \ref{T:raising}). In the case $L = S^{0}$ it has been proved in \cite{morita}.

\begin{pro}\label{P:1}
The following conditions are equivalent for any metrizable compactum $X$:
\begin{itemize}
\item[(i)]
$\dim_{L}Y \leq n$.
\item[(ii)]
There exists a map $f \colon X \to Y$ of a compactum $X$ with $\dim_{L}X = 0$ onto $Y$ such that $|f^{-1}(y)| \leq n+1$ for each $y \in Y$.
\end{itemize}
\end{pro}
\begin{proof}
(i) $\Longrightarrow$ (ii). Since $\dim_{L}Y \leq n$, there exists, by Theorem \ref{T:charact}, a map $p \colon Y \to I^{n}$ such that $\dim_{L}p = 0$. Let $\widetilde{Y} = p(Y)$. Since $\dim\widetilde{Y} \leq n$, there exist a zero-dimensional compactum $\widetilde{X}$ and a map $q \colon \widetilde{X} \to \widetilde{Y}$ onto $\widetilde{Y}$ such that $|q^{-1}(\widetilde{y})| \leq n+1$ for each $\widetilde{y} \in \widetilde{Y}$. Now let $X = \{ (\widetilde{x}, y) \in \widetilde{X} \times Y \colon q(\widetilde{x}) = p(y) \}$. Let also $f = \pi_{Y}|X \colon X \to Y$ and $g = \pi_{\widetilde{X}}|X \colon X \to \widetilde{X}$, where $\pi_{Y} \colon \widetilde{X} \times Y \to Y$ and $\pi_{\widetilde{X}} \colon \widetilde{X} \times Y \to \widetilde{X}$ denote the corresponding projections. In other words the following diagram

\[
\begin{CD}
X @> f >> Y\\
@V g VV @VV p V\\
\widetilde{X} @> q >> \widetilde{Y}\\
\end{CD}
\]

\noindent is a pullback square. Clearly fibers of the map $g$ are homeomorphic to the fibers of the map $p$ and therefore $\dim_{L}g = 0$. Since $\dim\widetilde{X} = 0$, we conclude, by Theorem
\ref {T:chihur}, that $\dim_{L}X = 0$. It only remains to note that the fibers of the map $f$ are homeomorphic to the fibers of the map $q$. Consequently, $|f^{-1}(y)| \leq n+1$ for each $y \in Y$.

(ii) $\Longrightarrow$ (i). This implication, as mentioned above, coincides with Theorem \ref{T:raising}.
\end{proof}


\subsection{Extensional properties of coronas}\label{SS:corona}
In this section we investigate dimensional properties of the Stone-\v{C}ech increment $\beta X \setminus X$ of a space $X$. 

We say that $L$ is an absolute extensor for a space $X$ with respect to the class of compact spaces if any map $f \colon A \to L$, defined on a closed subset $A \subseteq X$, admits an extension $\widetilde{f} \colon \operatorname{Cl}_{X}G \to L$, where $G$ is an open neighbourhood of $A$ in $X$ such that $X \setminus G$ is compact. In such a case we write $L \in {\rm AE}^{\rm c}(X)$.

\begin{thm}\label{T:corona}
If $L$ is finitely dominated, then following conditions are equivalent for any metrizable locally compact space $X$:
\begin{itemize}
\item[(a)]
$L \in {\rm AE}(\beta X \setminus X)$;
\item[(b)]
$L \in {\rm AE}^{\rm c}(X)$.
\end{itemize}
\end{thm}
\begin{proof}
(a) $\Longrightarrow$ (b). Let $f \colon A \to L$ be a map defined on a closed subset $A$ of $X$. Obviously, $\operatorname{Cl}_{\beta X}A$ is the Stone-\v{C}ech compactification $\beta A$ of $A$. Consequently, since $L$ is finitely dominated, there exists a map $g \colon \operatorname{Cl}_{\beta X}A \to L$ such that $g|A \simeq f$. Since $L \in {\rm AE}(\beta X \setminus X)$ there exists a map $g^{\prime} \colon (\beta X \setminus X) \to L$ such that $g^{\prime}|(\operatorname{Cl}_{\beta X}A \setminus X) = g|(\operatorname{Cl}_{\beta X}A \setminus X)$. Since $X$ is locally compact, the union $A \cup (\beta X \setminus X)$ is closed in $\beta X$. Therefore the map $g^{\prime\prime} \colon A \cup (\beta X \setminus X) \to L$, defined by letting
\[
g^{\prime\prime}(x) = 
\begin{cases}
g(x), \;\text{if}\; x \in A,\\
g^{\prime}(x),\;\text{if}\; x \in \beta X \setminus X ,
\end{cases}
\]

\noindent admits an extension $\widetilde{g} \colon V \to L$ onto an open set $V \subseteq \beta X$ such that $A \cup (\beta X \setminus X) \subseteq V$. Let $\widetilde{G}$ be an open subset of $\beta X$ such that $A \cup (\beta X \setminus X) \subseteq \widetilde{G} \subseteq \operatorname{Cl}_{\beta X}\widetilde{G} \subseteq V$. Let $G = \widetilde{G} \cap X$. Obviously, $A \subseteq G$, $X \setminus G = \beta X \setminus \widetilde{G}$ is compact and $\widetilde{g}|A = g|A \simeq f$. By the Homotopy Extension Theorem, $f$ admits an extension $\widetilde{f} \colon \operatorname{Cl}_{X}G \to L$ which proves that $L \in {\rm AE}^{\rm c}(X)$.

(b) $\Longrightarrow$ (a). Let $f \colon A \to L$ be a map, defined on a closed subset $A \subseteq \beta X \setminus X$. Let also $U$ and $V$ be open subsets of $\beta X$ such that $A \subseteq V \subseteq \operatorname{Cl}_{\beta X}V \subseteq U$ and $f$ admits an extension $f^{\prime} \colon U \to L$. The set $X \cap \operatorname{Cl}_{\beta X}V$ is non-empty and closed in $X$. Since $L \in {\rm AE}^{\rm c}(X)$, there exist an open set $G \subseteq X$ and a map $f^{\prime\prime}\colon \operatorname{Cl}_{X}G \to L$ such that $X \cap \operatorname{Cl}_{\beta X}V \subseteq G$, $X \setminus G$ is compact and $f^{\prime\prime}|(X \cap \operatorname{Cl}_{\beta X}V) = f^{\prime}|(X \cap \operatorname{Cl}_{\beta X}V)$. Compactness of $X \setminus G$ guarantees that $\operatorname{Cl}_{\beta X}G = \operatorname{Cl}_{\beta X}(\operatorname{Cl}_{X} G)$ is the Stone-\v{C}ech compactification of $\operatorname{Cl}_{X}G$ and $\beta X \setminus X \subseteq \operatorname{Cl}_{\beta X}G$. Since $\operatorname{Cl}_{X}G \subseteq \operatorname{Cl}_{X}G \cup \operatorname{Cl}_{\beta X}V \subseteq \operatorname{Cl}_{\beta X}G$ it follows that $\operatorname{Cl}_{\beta X}G$ is the Stone-\v{C}ech compactification of the sum $\operatorname{Cl}_{X}G \cup \operatorname{Cl}_{\beta X}V$ as well. Next consider the map $g \colon \operatorname{Cl}_{X}G \cup \operatorname{Cl}_{\beta X}V \to L$, defined by letting

\[
g(x) =
\begin{cases}
f^{\prime}(x), \;\text{if}\; x \in \operatorname{Cl}_{\beta X}V,\\
f^{\prime\prime}(x),\; \text{if}\; x \in \operatorname{Cl}_{X}G .
\end{cases}
\]

\noindent Since $L$ is finitely dominated, we can find a map $\widetilde{g} \colon \operatorname{Cl}_{\beta X}G \to L$ such that $\widetilde{g}|(\operatorname{Cl}_{X}G \cup \operatorname{Cl}_{\beta X}V) \simeq g|(\operatorname{Cl}_{X}G \cup \operatorname{Cl}_{\beta X}V)$. In particular, $\widetilde{g}|A \simeq f$. According to the Homotopy Extension Theorem (recall that $A \subseteq \beta X \setminus X \subseteq \operatorname{cl}_{\beta X}G$) the map $f$ admits an extension $\widetilde{f} \colon \beta X \setminus X \to L$. This proves that $L \in {\rm AE}(\beta X \setminus X)$.
\end{proof}

\section{Axiomatic characterization of $\dim_{L}$}\label{S:axiomatic}

A compact metrizable space $X$ such that $\operatorname{dim}_{L}X = n$, is a generalized Cantor $L_{n}$-manifold if there is no closed subset $Y$ of $X$, satisfying the inequality $\dim_{L} Y \leq n-2$, such that $X \setminus Y$ is disconnected.

\begin{lem}\label{L:196}
Let $f,g \colon X \to \Sigma^{n}L$ be continuous maps of a separable metrizable space $X$. If  $\dim_{L}\left(\{ x \in X \colon f(x) \neq g(x) \}\right) \leq n-1$, then $f \simeq g$.
\end{lem}
\begin{proof}
Let $Y = \{ x \in X \colon f(x) \neq g(x)\}$ and consider the map 

\[ h \colon  \left( X \times \{ 0,1\}\right) \cup  \left( (X \setminus Y) \times [0,1]\right) \to \Sigma^{n}L , \]

\noindent defined by letting

\[ h(x,t) = 
\begin{cases}
f(x) , \text{if}\; (x,t) \in X \times \{ 0\} ,\\
f(x) , \text{if}\; (x,t) \in \left( X \setminus Y\right) \times [0,1] ,\\
g(x) , \text{if}\; (x,t) \in X \times \{ 1 \} .\\
\end{cases}
\]

Note that $X \times [0,1] \setminus \left( (X\times \{ 0,1\} ) \cup (X\setminus Y)\times [0,1]\right) \subseteq  Y \times [0,1]$. By our assumption, $\dim_{L}Y \leq n-1$. In other words, $\Sigma^{n-1}L \in {\rm AE}(Y)$. By Theorem \ref{T:suspension}, $\Sigma^{n}L = \Sigma\left(\Sigma^{n-1}L\right) \in {\rm AE}(Y \times [0,1])$. Obviously this suffices to conclude that the map $h$ admits an extension $H \colon X \times [0,1] \to \Sigma^{n}L$ which provides a needed homotopy between the maps $f$ and $g$.
\noindent
\end{proof} 

\begin{thm}\label{T:manifold}
Let $X$ be a metrizable compactum. If $\dim_{L}X = n \geq 2$, then $X$ contains a generalized Cantor $L_{n}$-manifold.
\end{thm}
\begin{proof}
Since $\dim_{L}X = n$ it follows that $\Sigma^{n}L \in {\rm AE}(X)$, but $\Sigma^{n-1}L \notin {\rm AE}(X)$. Thus there exists a map $f \colon F \to \Sigma^{n-1}L$, defined on a closed subset $F$ of $X$, which is not extendable over $X$. Let ${\mathcal F}$ denote the partially ordered set (by inclusion) of all closed subsets $Y \subseteq X$ such that the map $f$ cannot be extended over $F \cup Y$. This set is non-empty , since $X \in {\mathcal F}$. Using the compactness of $X$, the fact that $\Sigma^{n-1}L$ is an ${\rm ANR}$-space and the Kuratowski-Zorn lemma, we conclude, following the standard argument, that ${\mathcal F}$ contains a maximal element $Y$. In other words, there exists a closed subset $Y \subseteq X$ such that $f$ is not extendable over $F \cup Y$, but is extendable over $F \cup Y^{\prime}$ for any proper closed subset $Y^{\prime}$ of $Y$. We claim that $Y$ is a generalized Cantor $L_{n}$-manifold. First note that $\dim_{L}Y \leq \dim_{L}X = n$. Next suppose that $Y$ is represented as the union of its two proper closed subsets $Y_{1}$ and $Y_{2}$. The proof will be completed if we show that $\dim_{L}(Y_{1}\cap Y_{2}) \geq n-1$. Assume the contrary, i.e. that $\dim_{L}(Y_{1}\cap Y_{2}) \leq n-2$. By construction, $f$ can be extended over $f_{k} \colon F \cup Y_{k} \to \Sigma^{n-1}L$, for each $k = 1,2$. Since $\{ x \in F \cup (Y_{1}\cap Y_{2}) \colon f_{1}(x) \neq f_{2}(x) \} \subseteq Y_{1}\cap Y_{2}$, we conclude, by Lemma \ref{L:196} and by our assumption, that $f_{1}|(F \cup(Y_{1}\cap Y_{2})) \simeq f_{2}|(F\cup (Y_{1} \cap Y_{2}))$. The Homotopy Extension Theorem guarantees that $f_{1}|(F \cup(Y_{1}\cap Y_{2}))$ admits an extension $\widetilde{f}_{1} \colon F \cup Y_{2} \to \Sigma^{n-1}L$. Then the map $g \colon F \cup Y \to \Sigma^{n-1}L$, which coincides with $f_{1}$ on $F \cup Y_{1}$ and with $\widetilde{f}_{1} \colon F \cup Y_{2}$, is continuous and extends $f$. This contradicts the choice of $f$ and $Y$ and completes the proof.
\end{proof}


Below let ${\mathcal K}$ denote the class of finite dimensional in the sense of $\dim$ metrizable compact spaces. Similarly ${\mathcal K}_{L}$ denotes the class of finite dimensional in the sense of $\dim_{L}$ metrizable compacta.

Let $d \colon {\mathcal K} \to \{ -1, 0,1,2,\dots \}$ be an integer-valued function which assigns same values to any pair of homeomorphic spaces.
In 1932 Alexandrov gave  (see \cite[Chapter 5, $\S$ 10, Theorem 19]{alpas}
the following characterization of the dimension function $\dim$.

{\bf Alexandrov's Axiomatization.} {\em The Lebesgue covering dimension $\dim$ is the only function $d$ which satisfies conditions {\rm (A1)--A(4)} in the class ${\mathcal K}$:
\begin{description}
\item[(A1) -- normalization axiom]
$d(\emptyset ) = -1$, $d(I^{n}) = n$ for $n = 0,1,2,\dots$.
\item[(A2) -- sum axiom]
If the space $X \in {\mathcal K}$ is represented as the union of two closed subspaces $X_{1}$ and $X_{2}$, then $d(X) = \max \{ d(X_{1}), d(X_{2})\}$.
\item[(A3) -- Poincar\'{e}'s axiom]
For every $X \in {\mathcal K}$  with $|X| >1$ there exists a closed set $X^{\prime} \subseteq X$ separating $X$ and such that $d(X^{\prime}) < d(X)$.
\item[(A4) -- Brouwer's axiom]
For every space $X \in {\mathcal K}$ there exists an open cover $\omega$ such that if $f \colon X \to Y$ is an $\omega$-map of $X$ onto a space $Y \in {\mathcal K}$, then $d(X) \leq d(Y)$.
\end{description}}

It is not hard to see (see, for instance, \cite[footnote on p. 976]{scepin}) that the Brouwer's axiom can be replaced  by either of the following conditions:
\begin{description}
\item[(A5) -- continuity axiom]
If ${\mathcal S} = \{ X_{k}, p_{k}^{k+1}\}$ is an inverse sequence of spaces from ${\mathcal K}$, then $d(\lim{\mathcal S}) \leq \sup\{ d(X_{k})\}$.
\item[(A6) -- Hurewicz's axiom]
If there exists a map $f \colon X \to I^{n}$ such that $d(f^{-1}(y)) =0$ for every $y \in f(X)$, then $d(X) \leq n$.
\end{description}

 In order to characterize extraordinary dimension function $\dim_{L}$ with classifying complex $L$ we need to adjust some of the above axioms and replace Brouwer's axiom by Hurewicz's axiom.

\begin{thm}\label{T:axiomatics}
The dimension $\dim_{L}$ is the only function,
defined on the class ${\mathcal K}_{L}$, which satisfies the following axioms:
\begin{description}
\item[(C1)-- normalization axiom]
$d(X) \in \{ 0,1,2,\dots \}$ and $d(X) = 0$ if and only if $L \in {\rm AE}(X)$. 
\item[(C2) -- monotonicity axiom]
If $A$ is a closed subspace of $X$, then $d(A) \leq d(X)$.
\item[(C3) -- Poincar\'{e}'s axiom]
If $d(X) > 0$, then there exists a closed subspace $A$ in $X$ separating $X$ and such that $d(A) < d(X)$.
\item[(C4) -- Hurewicz's axiom]
If there exists a map $f \colon X \to I^{n}$ such that $d(f^{-1}(y)) =0$ for every $y \in f(X)$, then $d(X) \leq n$.
\end{description} 
\end{thm}
\begin{proof}
First let us show that $d(X) \leq \dim_{L}X$ for every $X \in {\mathcal K}_{L}$. If $\dim_{L}X \leq n$, then, by Proposition  \ref{P:cubes}, there exists a map $f \colon X \to I^{n}$ such that $\dim_{L}f^{-1}(y) = 0$ for each $y \in f(X)$. By (C1), $d(f^{-1}(y) = 0$ for each $y \in f(X)$. Consequently, by (C4), $d(X) \leq n$.

Next we show that $\dim_{L}X \leq d(X)$ for every $X \in {\mathcal K}_{L}$. If $d(X) = 0$, then, by (C1), $\dim_{L}X = 0$. Suppose that the inequality $\dim_{L}Y \leq d(Y)$ has been proved for spaces $Y \in {\mathcal K}_{L}$ with $d(Y) \leq n-1$, $n\geq 1$, and consider a space $X \in {\mathcal K}_{L}$ such that $d(X) = n$. Next assume the contrary, i.e. $\dim_{L}X = m > n$. Note that $m \geq 2$. By Theorem \ref{T:manifold}, $X$ contains a Cantor $L_{m}$-manifold $Z$. By (C2), $d(Z) \leq d(X) = n$. Note that $d(Z) >0$ (otherwise, by (C1), we get $0 = \dim_{L}Z = \dim_{L}X = m \geq 2$). By (C3), there exists a closed subspace $Y$ of $Z$ which separates $Z$ and such that $d(Y) < d(Z) \leq n$. Then $d(Y) \leq n-1$. By the inductive hypothesis, $\dim_{L}Y \leq d(Y) \leq n-1 < m-1$. Consequently the closed subspace $Y$ of dimension $\dim_{L}Y < m-1$ separates the  Cantor $L_{m}$-manifold $Z$. This contradiction completes the proof.
\end{proof}


\section{Appendix: Spectral characterizations of the relation $L \in {\rm AE}(X)$}\label{S:app}

In this section we present spectral characterizations of the relation $L \in {\rm AE}(X)$ which have been used in the proofs throughout this paper. Definitions of concepts related to inverse spectra can be found in \cite{chibook}.

The standard situation we would like to analyze is as follows. We are given a realcompact and $z$-embedded subspace $Y$ of a realcompact space $X$. Also we have a Polish spectrum ${\mathcal S}_{X} = \{ X_{\alpha}, p_{\alpha}^{\beta}, A\}$ such that $X = \lim{\mathcal S}_{X}$. Let us see how the relations $L \in {\rm AE}(X)$ and $L \in {\rm AE}(Y)$ can be characterized in term of the given spectrum ${\mathcal S}_{X}$. Answer to the first question has been given in the following statement.

\begin{thm}[\cite{chi3}, Theorem 4.4]\label{T:unc}
Let $L$ be a Polish ${\rm ANR}$-space. Then the following conditions are equivalent for any realcompact space $X$ and Polish spectrum ${\mathcal S}_{X} = \{ X_{\alpha}, p_{\alpha}^{\beta}, A\}$ with $X = \lim{\mathcal S}_{X}$:
\begin{itemize}
\item[(a)]
$L \in {\rm AE}(X)$.
\item[(b)]
There exists a cofinal and $\omega$-complete subset $A^{\prime}$ of the indexing set $A$ such that $L \in {\rm AE}(X_{\alpha})$ for each $\alpha \in B$.
\end{itemize}
\end{thm}

Let us now analyze the relation $L \in {\rm AE}(Y)$. Of course, since $Y$ itself is a realcompact space, one can apply Theorem \ref{T:unc} to any Polish spectrum ${\mathcal S}_{Y} = \{ Y_{\alpha}, q_{\alpha}^{\beta}, A\}$ with $Y = \lim{\mathcal S}_{Y}$ and find a cofinal and $\omega$-complete subset $A^{\prime\prime} \subseteq A$ such that $L \in {\rm AE}(Y_{\alpha})$ for each $\alpha \in A^{\prime\prime}$. Problem is that the spectrum ${\mathcal S}_{Y}$ in no way reflects the fact that $Y$ is a subspace of $X$: we cannot assume, even if $A^{\prime} = A^{\prime\prime}$, that $Y_{\alpha}$ is a subspace of $X_{\alpha}$ for sufficiently many indices. A logical way to establish such a connection would be to consider the induced spectrum ${\mathcal S} = \{ p_{\alpha}(Y), p_{\alpha}^{\beta}|p_{\beta}(Y), A\}$ and then to apply the Spectral Theorem \cite[Theorem 1.3.6]{chibook} to the spectra ${\mathcal S}_{Y}$, ${\mathcal S}$ and to the inclusion map $i \colon Y \hookrightarrow X$.
In order to be able to proceed this way we need to know that the spectrum ${\mathcal S}$ is also factorizing and $\omega$-continuous. While it is indeed factorizing (this follows from \cite[Propositions 1.1.22 and 1.1.24]{chibook}), it is clearly not $\omega$-continuous. Nevertheless it is still possible to extract some information about the relation $L \in {\rm AE}(Y)$ from the spectrum ${\mathcal S}$. We record this information in the following statement proof of which can be extracted from the proofs of \cite[Theorem 1.3.6]{chibook} and \cite[Theorems 4.4 and 6.5]{chi3}.

\begin{pro}\label{P:nonpolish}
Let $Y$ be a realcompact and $z$-embedded subspace of a realcompact space $X$. Let also ${\mathcal S}_{X} = \{ X_{\alpha}, p_{\alpha}^{\beta}, a\}$ be a Polish spectrum such that $X = \lim{\mathcal S}_{X}$ and $L$ be a Polish ${\rm ANR}$-space. Then
\begin{itemize}
\item[(a)]
If $L \in {\rm AE}(Y)$, then there exists a cofinal and $\omega$-complete subset $A^{\prime} \subseteq A$ such that $L \in {\rm AE}(p_{\alpha}(Y))$.
\item[(b)]
If $L \in {\rm AE}(p_{\alpha}(Y))$ for each $\alpha \in A^{\prime}$, where $A^{\prime}$ is a cofinal and $\omega$-complete subset of $A$, then $L \in {\rm AE}(Y)$.
\end{itemize}
\end{pro}

The following statement for metrizable spaces appears in Theorem \ref{T:menger}. 

\begin{thm}\label{T:join}
Let $X$ and $Y$ be $z$-embedded subspaces of the union $Z = X \cup Y$. If $L \in {\rm AE}(X)$ and $K \in {\rm AE}(Y)$, then $L\ast K \in {\rm AE}(Z)$.
\end{thm}
\begin{proof}
Let $\upsilon Z$ be the Hewitt realcompactification of $Z$ and let $\widetilde{X}$ denote the intersection of all functionally open subsets of $\upsilon Z$, containing $X$. Note that $\widetilde{X}$ is homeomorphic to $\upsilon X$ and is a $z$-embedded in $\upsilon Z$ (see \cite[Proposition 1.1.24]{chibook}). Let $\widetilde{Y}$ has the similar meaning. By Theorem \ref{T:hewitt}, $L \in {\rm AE}(\widetilde{X})$ and $K \in {\rm AE}(\widetilde{Y})$. Note also that $\upsilon Z = \widetilde{X} \cup \widetilde{Y}$. Next consider any Polish spectrum ${\mathcal S}_{\upsilon Z} = \{ Z_{\alpha}, p_{\alpha}^{\beta}, A\}$ such that $\upsilon Z = \lim{\mathcal S}_{\upsilon Z}$. By Proposition \ref{P:nonpolish}(a), there exist cofinal and $\omega$-complete subsets $A_{X}, A_{Y} \subseteq A$ such that $L \in {\rm AE}(p_{\alpha}(\widetilde{X}))$ for each $\alpha \in A_{X}$ and $K \in {\rm AE}(p_{\alpha}(\widetilde{Y}))$ for each $\alpha \in A_{Y}$. By \cite[Proposition 1.1.27]{chibook}, the intersection $A_{X} \cap A_{Y}$ is cofinal and $\omega$-closed in $A$. Since $Z_{\alpha}$ is metrizable (and even separable), we conclude, by Theorem \ref{T:menger}, that $L \ast K \in {\rm AE}(p_{\alpha}(\upsilon Z))$ for each $\alpha \in A_{X} \cap A_{Y}$  (note here that $p_{\alpha}(\upsilon Z)$ generally speaking is a proper subset of $Z_{\alpha}$ and consequently we are not able to conclude that $L \ast K \in {\rm AE}(Z_{\alpha})$). Proposition \ref{P:nonpolish}(b) guarantees that $L\ast K \in {\rm AE}(\upsilon Z)$. Once again applying Theorem \ref{T:hewitt} we conclude that $L\ast K \in {\rm AE}(Z)$.
\end{proof}

Similar considerations prove the following corollary.

\begin{thm}[\cite{chi3}, Proposition 6.8]\label{T:monotone}
Let $Y$ be a $z$-embedded subset of a space $X$. If $L \in {\rm AE}(X)$, then $L \in {\rm AE}(Y)$. In particular, the latter holds if:
\begin{itemize}
\item[(a)]
$Y$ is an $F_{\sigma}$-subset of a normal space $X$;
\item[(b)]
$Y$ is any subset of a perfectly normal space $X$;
\item[(c)]
$Y$ is an open or a dense subset of a perfectly $\kappa$-normal space $X$.
\item[(d)]
$Y$ is a Lindel\"{o}f subspace of a completely regular and Hausdorff space.
\end{itemize}
\end{thm}

\end{document}